\documentclass{article}

\usepackage[american]{babel}
\usepackage[utf8]{inputenc}

\usepackage{amsmath, amssymb}   
\usepackage{amsthm}
\usepackage{mathrsfs}
\usepackage{mathtools}
\usepackage{stix}
\usepackage{color}
\usepackage{changes}

\usepackage{makeidx}
\usepackage{acronym}
\usepackage{algorithm}
\usepackage{algpseudocode}
\usepackage{longtable} 

\usepackage{enumitem}
\usepackage{adjustbox} 
\usepackage{tikz}
\usetikzlibrary{calc, shapes, positioning, decorations}

\usepackage{pgfplots}
\usepackage{pgfplotstable}
\pgfplotsset{compat=1.9}

\usepackage[T1]{fontenc}
\usepackage{titlesec}

\theoremstyle{plain}
\newtheorem{theorem}{Theorem}[section]
\newtheorem{corollary}[theorem]{Corollary}
\newtheorem{proposition}[theorem]{Proposition}
\newtheorem{lemma}[theorem]{Lemma}

\newtheorem{example}[theorem]{Example}

\theoremstyle{definition}
\newtheorem{definition}[theorem]{Definition}

\newtheorem{remark}[theorem]{Remark}

\newcommand{\N}{\mathbb{N}}
\newcommand{\C}{\mathbb{C}}
\newcommand{\Z}{\mathbb{Z}}

\newcommand{\Q}{\mathbb{Q}}

\newcommand{\F}{\mathbb{F}}

\newcommand{\diff}[2]{#1 \cdot #2^{-1}}

\renewcommand{\paragraph}[1]{\leavevmode\vspace{0,2cm}\newline{\it #1}\\

}
\renewcommand{\tilde}{\widetilde}
\let\oldReturn\Return
\renewcommand{\Return}{\State\oldReturn}

\DeclareMathOperator{\wt}{wt}

\DeclareMathOperator{\ord}{ord} 
\DeclareMathOperator{\TITO}{TITO}
\DeclareMathOperator{\Tr}{Tr} 

\DeclareMathOperator{\Aut}{Aut} 

\DeclareMathOperator{\id}{id}

\begin{document}

\title{Formal self duality}

\author{Lukas Kölsch\thanks{lukas.koelsch@uni-rostock.de, Institute of Mathematics, University of Rostock, 18051 Rostock, Germany}, Robert Schüler\thanks{robert.schueler2@uni-rostock.de, Institute of Mathematics, University of Rostock, 18051 Rostock, Germany}}

\maketitle

\begin{abstract}
We study the notion of formal self duality in finite abelian groups. Formal duality in finite abelian groups has been proposed by Cohn, Kumar, Reiher and Schürmann. In this paper we give a precise definition of formally self dual sets and discuss results from the literature in this perspective. Also, we discuss the connection to formally dual codes. We prove that formally self dual sets can be reduced to primitive formally self dual sets similar to a previously known result on general formally dual sets. Furthermore, we describe several properties of formally self dual sets. Also, some new examples of formally self dual sets are presented within this paper. Lastly, we study formally self dual sets of the form $\{(x,F(x)) \ : \ x\in\F_{2^n}\}$ where $F$ is a vectorial Boolean function mapping $\F_{2^n}$ to $\F_{2^n}$.
\end{abstract}

\section{Introduction}

In this paper we study formal duality as introduced by Cohn, Kumar, Reiher and Schürmann in \cite{cohn2009ground}.
For an outline of the relations of this concept to formal duality of codes, we refer to Section \ref{sec:connection_to_codes}.
Formal duality has been introduced 
in relation to energy minimization problems and has subsequently been studied using finite abelian groups in \cite{cohn2014formal}.

It can be defined as follows:

\begin{definition}\label{def:formal-duality}
Let $G$ be some (multiplicative) finite abelian group and
$\hat G$ be its dual group, i.e., the group of homomorphisms from $G$ to $\C^\ast$. Two sets $S \subset G$ and $T \subset \hat G$
form  a formally dual pair, if for all
$\chi\in\hat G$ we have

\begin{equation}
\left|\chi(S)\right|^2 = \frac{|S|^2}{|T|}\nu_T(\chi),\label{eq:def_fd}
\end{equation}
or equivalently if for all $g\in G$ we have
\begin{equation}
\left|g(T)\right|^2 = \frac{|T|^2}{|S|}\nu_S(g),\label{eq:equivalent_def_fd}
\end{equation}
where 
$$\nu_T(\chi) = |\{(\phi,\psi)\in T\times T \ : \ \diff \phi\psi = \chi\}|$$ is called the \emph{weight enumerator} of $T$ and we use the notation $\chi(S) = \sum_{x\in S} \chi(x)$ and $g(T) = \sum_{\chi\in T} \chi(g)$.
A set $S$ is called a \emph{formally dual set} if there is a set $T$ such that $S$ and $T$ form a formally dual pair.
\end{definition}

Note that $\hat{\hat G}$ is canonically isomorphic to $G$ and $S$ and $T$ form a formally dual pair (in $G$ and $\hat G$) if and only if $T$ and $S$ form a formally dual pair (in $\hat G$ and $\hat{\hat G}$). Furthermore, formal duality can be seen as a generalization of relative difference sets, i.e. sets $S\subset G$ such that there is a subgroup $N\leq G$ and an integer $\lambda$ with $\nu_S(g) = |G|$ if $g = 1$, $\nu_S(g) =  0$ if $g\in N$ and $\nu_S(g) = \lambda$ otherwise.

The characterization of such formally dual sets is an interesting open question which several authors studied before.
More information about formal duality in cyclic groups can be found in \cite{schuler2017Cyclic}, \cite{xia2016classification}, \cite{malikiosis2017cyclic}. A comprehensive analysis of formal duality in general abelian groups is given in \cite{li2018abelian}, \cite{schuler2019phd} and some more examples are discussed in \cite{li2019constructions}, \cite{li2020direct_construction}.

In order to define formal self duality, we have to regard the set $T$ of Definition~\ref{def:formal-duality} as a subset of $G$.
A finite abelian group is always isomorphic to its dual group. Thus, by choosing an isomorphism $\Delta: G\rightarrow \hat G$ we can get rid of the dual group in Definition \ref{def:formal-duality}. Equivalently we can choose a \emph{pairing}, that is a nondegenerate bilinear form $\left<\cdot,\cdot\right>:G\times G \rightarrow \C^\ast$. Note that any isomorphism $\Delta$ defines a pairing $\left<\cdot,\cdot\right>_\Delta$ by $\left<a,b\right>_\Delta = [\Delta a](b)$ and any pairing $\left<\cdot,\cdot\right>$ defines an isomorphism by $x\mapsto \chi$ such that $\chi(y) = \left<x,y\right>$ for all $y\in G$.

Using this notion we can define formal duality under isomorphisms and formal self duality:

\begin{definition}[{\cite[Definition 2.7]{li2018abelian}}]\label{def:fd_under_isomorphism}
Let $G$ be a finite abelian group and \mbox{$\Delta: G\rightarrow \hat G$} be an isomorphism. Then $S\subset G$ and $T\subset G$ form a \emph{formally dual pair under} $\Delta$ if $S$ and $\Delta(T)$ form a formally dual pair (in $G$ and $\hat G$). Alternatively, we say $S$ and $T$ form a formally dual pair under the pairing $\left<\cdot,\cdot\right>_\Delta$. We call a set $S$ \emph{formally self dual}, when $S$ is formally dual to itself under some isomorphism.
\end{definition}

Formal self duality has been briefly studied in \cite{xia2016classification} and some examples are given in \cite{cohn2014formal}, \cite{li2018abelian}. 

In this paper we continue the study of formal self duality.
An important tool in the study of formal duality is the reduction to so called primitive formally dual sets (see Theorem \ref{thm:reduce_to_primitive}). We prove an analog result for formal self duality. Furthermore, we state some equivalent formulations of formal self duality using the even set theory introduced in \cite{li2018abelian}.
Moreover, we state four new examples of primitive formally dual sets in groups of order $64$ that happen to be formally self dual. Note that in the set of all groups with order no more than $63$ we are able to characterize formally dual sets (see \cite[Table A.1]{li2018abelian}). Thus, these new examples are a step towards increasing the smallest group size with incomplete results.

The paper is organized as follows:
We start with a comparison of formal duality in finite abelian groups and formal duality of codes in Section \ref{sec:connection_to_codes}.
In Section \ref{sec:preliminaries} we state previously known facts needed to follow the rest of the paper. We also give an overview of examples of formal self duality known from the literature. In Section \ref{sec:results} we state and prove the results about formal self duality mentioned before and Section \ref{sec:examples} contains the new examples. In Section \ref{sec:boolean_functions}, we discuss the relations of formal duality to Boolean vectorial functions and study formally self dual sets of the form $\{(x,F(x)) \ : \ x\in\F_{2^n}\}$. We conclude this paper by some open questions in Section \ref{sec:open_questions}.

\section{Connection to formal duality of codes}\label{sec:connection_to_codes}

In this section we point out an interesting connection between formal duality of sets in finite abelian groups as defined in Definition~\ref{def:formal-duality} and the notion of \emph{formal dual codes}.
The possible connection between these to concepts has been brought to our attention by Claude Carlet.
Let $C \subseteq \F_q^n$ be a code with $q=p^k$. We denote the weight of $c \in C$  (i.e. the number of non-zero entries) by $\wt(c)$ and the Hamming distance of $c_1,c_2$ by $d_H(c_1,c_2)=\wt(c_1-c_2)$. Further,  we define the weight enumerator $W_C$ and distance enumerator $D_C$ by
\[W_C(X,Y) = \sum_{c \in C}X^{n-\wt(c)}Y^{\wt(c)} \text{ and } D_C(X,Y) = \frac{1}{|C|}\sum_{c_1,c_2 \in C}X^{n-d_H(c_1,c_2)}Y^{d_H(c_1,c_2)}.\]
If $C$ is a linear code, it is elementary that $W_C=D_C$. A key result in the theory of linear codes is the MacWilliams identity (see e.g.~\cite{coding}), which gives a relation between the weight/distance enumerator of a linear code $C$ and its dual $C^\perp = \{v \in \F_q^n \colon v \cdot c = 0 \text{ for all } c \in C\}$, where the multiplication $\cdot$ denotes the usual dot product on $\F_q^n$. The MacWilliams identity then states
\[\frac{1}{|C|}W_C(X+(q-1)Y,X-Y) = W_{C^\perp}(X,Y),\]
and, of course, the same relation equivalently for the distance enumerator $D_C$. Extending the notion of dual codes, we say that two non-linear codes $C,C' \subseteq \F_q^n$ that nevertheless satisfy the MacWilliams identity in the sense that  
\[\frac{1}{|C|}W_C(X+(q-1)Y,X-Y) = W_{C'}(X,Y)\]
have dual weight enumerators. Equivalently, if the same relation holds for the distance enumerators, we say that $C,C'$ have dual distance enumerators. If $C,C'$ have both dual weight and distance enumerators, we call them \emph{formal dual codes}. In the case of $C=C'$, we can also speak of formal self dual codes. Probably the most famous case of formal dual codes are the binary Kerdock and Preparata codes~\cite[Chapter 15]{coding}. The breakthrough result in~\cite[Theorems 2,3]{z4linearity} shows in fact that binary codes constructed via the Gray mapping from $\Z_4$-linear codes and their $\Z_4$-duals are always formal dual codes in the binary setting. For a more detailed treatment of this connection between $\Z_4$-linear codes and formal dual codes, we also refer to~\cite{carlet_formalduality}. There is also a connection between the notion of formal duality for codes and the notion of formal duality we introduced for finite abelian groups in Definition~\ref{def:formal-duality}: Following the steps of the classical proof of the MacWilliams identity~\cite[Theorem 13]{coding}, we have for arbitrary codes $C \subseteq \F_q^n$
\begin{equation}
D_C(X+(q-1)Y,X-Y) = \frac{1}{|C|}\sum_{x \in \F_q^n} \left|\sum_{c \in C} \chi_x(c)\right|^2 X^{n-\wt(x)}Y^{\wt(x)},
\label{eq:dist_enumerator}
\end{equation}
where $\chi_x(c) = \zeta_p^{\Tr( x\cdot c )}$, $\zeta_p=e^{2\pi i/p}$ is a primitive $p$-th root of unity and $\Tr \colon \F_q \rightarrow \F_p$ the absolute trace mapping $x \mapsto x+x^p+\dots+x^{p^{k-1}}$. Indeed, the MacWilliams identity for linear codes follows directly from Eq.~\eqref{eq:dist_enumerator} since in that case $\sum_{c \in C} \chi_x(c)=|C|$ if $x \in C^\perp$ and $0$ otherwise. Let us now assume that there is a set $C'\subseteq \F_q^n$ that is formal dual to $C$ as subsets of $(\F_q^n,+)$ in the sense of Definition~\ref{def:fd_under_isomorphism} under the canonical isomorphism $x \mapsto \chi_x$. Then, by Eq.~\eqref{eq:def_fd}, for each $x \in \F_q^n$ we have $\left|\sum_{c \in C}\chi_x(c)\right|^2 = \frac{|C|^2}{|C'|}\nu_{C'}(x)$. Then, using Eq.~\eqref{eq:dist_enumerator}, we have
\begin{align*}
D_C(X+(q-1)Y,X-Y) &= \frac{|C|}{|C'|}\sum_{x \in \F_q^n} \nu_{C'}(x) X^{n-\wt(x)}Y^{\wt(x)} \\
&= \frac{|C|}{|C'|}\sum_{c_1,c_2 \in C'}  X^{n-d_H(c_1,c_2)}Y^{d_H(c_1,c_2)} = |C| \cdot D_{C'}(X,Y).
\end{align*}
In particular, formal dual sets in $\F_q^n$ under the canonical isomorphism always yield codes that have dual distance enumerators.

It is however not the case, that formal dual sets in groups of the form $\F_q^n$ always yield formal dual codes.
Indeed, formal duality in abelian groups is invariant under translations. But formal dual codes $C,C'$ always have to contain $(0,\dots,0)$ (since $\frac{1}{|C|}W_C(X+(q-1)Y, X-Y)$ always contains the term $X^n$ and $W_C(X,Y)$ contains the term $X^n$ only if $C$ contains $(0,\dots,0)$). Thus, formal duality of codes is not invariant under translations. A simple counterexample can thus be constructed by translating any formal dual subset of a suitable group in such a way that it does not contain the zero vector. But there are also formal dual sets in groups of the form $\F_q^n$ which contain $(0,\dots,0)$ and still are not formal dual codes.
For example, consider the set 
\begin{align*}
	C = \{&(0, 0, 0, 0), (0, 1, 0, 1), (0, -1, 0, 1), (1, 0, 0, -1), (1, 1, -1, 0), (1, -1, 1, 0),\\
	&(-1, 0, 0, -1), (-1, 1, 1, 0), (-1, -1, -1, 0)\}\in\F_3^4.
\end{align*}
Note that $C$ forms with $C' = \Delta(C)$ a formally dual pair in $\F_3^4$ under the canonical isomorphism, where $\Delta(x_1,x_2,x_3,x_4) = (x_3,-x_4,x_1,-x_2)$. Further, $W_C(X,Y) = W_{C'}(X,Y)$ since the weight is clearly invariant under $\Delta$. Moreover, $(0,0,0,0)\in C\cap C'$.
But the MacWilliams identity does not hold since
\begin{align*}
	W_C(X,Y) &= X^4 + 4X^2Y^2 + 4XY^3\\
	\frac{1}{|C'|}W_{C'}(X+2Y, X-Y) &= X^4 + \frac 43 X^3Y + 4 XY^3 + \frac 83 Y^4.
\end{align*}
Therefore, $C$ and $C'$ are not formal dual codes.

On the other hand, it is also not the case that formal dual codes always yield formal dual sets.

As a counterexample, consider the two $\Z_4$-linear codes $C = \langle (2,1,3,1),(1,2,1,3) \rangle$, $C' = \langle (1,3,1,0),(3,1,0,1) \rangle$. It is easy to check that $C,C'$ are $\Z_4$-duals, i.e. $C' = \{x \in \Z_4^4 \colon x \cdot c = 0 \text{ for all } c \in C\}$, where $\cdot$ again denotes the usual dot product in $\Z_4^4$. Then, by \cite[Theorems 2,3]{z4linearity}, the associated binary codes $\phi(C),\phi(C') \subseteq \F_2^{8}$, where $\phi\colon \Z_4^4 \rightarrow \F_2^8$ is the Gray mapping (see \cite{z4linearity}), are formal dual codes. We explicitly state the resulting codes for convenience:
\begin{align*}
	\phi(C) =& \{(0,0,0,0,0,0,0,0),(1,0,1,0,1,1,0,1),(0,1,1,1,0,1,1,1),(1,1,0,1,1,0,1,0), \\
	&(0,1,0,1,1,1,1,0),(1,0,1,1,1,0,1,1),(1,1,1,0,0,1,0,1),(1,1,0,0,0,0,0,0), \\
	&(1,1,0,0,1,1,0,0),(0,0,0,0,1,1,0,0),(0,0,0,1,0,1,1,0),(0,1,1,0,0,0,0,1), \\
	&(0,1,1,1,1,0,1,1),(0,0,1,0,1,0,0,1),(1,0,0,1,0,0,1,0),(1,0,1,1,0,1,1,1)\}, \\
	\phi(C') =& \{(0,0,0,0,0,0,0,0),(0,1,0,0,1,0,1,0),(1,1,1,0,1,1,1,0),(1,0,1,0,0,1,0,0),\\
	&(1,0,0,0,0,1,0,1),(1,1,0,1,1,1,0,1),(0,1,0,1,1,0,0,0),(0,0,0,0,0,0,1,1), \\
	&(1,0,0,1,0,1,1,1),(1,1,0,1,1,1,1,0),(0,1,1,0,1,0,1,1),(0,0,1,1,0,0,1,1),\\
	&(1,0,1,1,0,1,1,0),(1,1,1,0,1,1,0,1),(0,1,1,1,1,0,0,1),(0,0,1,1,0,0,0,0)\}.
\end{align*}
Note that both codes are non-linear over $\F_2$. However, the condition for formally dual sets is not fulfilled by $\phi(C),\phi(C')$:
Indeed, as can be checked easily by computer, we have $\sum_{c \in \phi(C)}\chi(c) = 0$ for only $216$ choices of $\chi$ while $\nu_{\phi(C')}(x) = 0$ for $228$ choices of $x$, violating Eq.~\eqref{eq:def_fd} for any choice of isomorphism as in Definition \ref{def:fd_under_isomorphism}.

However, there are also sets which are simultaneously formally self dual sets and formal self dual codes. For example, when $p$ is an odd prime, the set $C = \{(x,x^2) \ : \ x\in\F_p\}$ is a formally self dual set under the pairing $\langle (a,b),(x,y)\rangle = \zeta_p^{\Tr(ay + bx)}$ in the additive group of $\F_p^2$ (see \cite[Theorem 3.2]{cohn2014formal}).
This set is also a formally self dual code: Indeed $\wt((x,x^2)) = \begin{cases} 0 \text{ if } x = 0,\\ 2 \text{ otherwise,}\end{cases}$ and therefore $C$ has a self dual weight enumerators as
$$\frac{1}{|C|}W_C(X+(p-1)Y, X-Y) = X^2 + (p-1)Y^2 = W_{C}(X,Y).$$
Similar to the discussion above, we can see that $C$ also has a self dual distance enumerator and is therefore a formally self dual code.

Of course, the motivation of the study of formal dual sets in finite abelian groups is quite different from the motivations in coding theory. For instance, the coding theoretical properties like minimum distance are one of the main reasons to investigate the formal dual Kerdock and Preparata codes, but are not of particular interest in our study.

\section{Preliminaries}\label{sec:preliminaries}

In this section we give a brief summary of needed results from the literature.
Recall, that the definition of formal self duality, i.e., Definition \ref{def:fd_under_isomorphism}, uses an isomorphism from $G$ to $\hat G$ or a pairing in $G$. Often, a very specific paring is used.
If a group is given in the form $G = \Z_{n_1}\times\Z_{n_2}\times\dots\times\Z_{n_m}$ we define the \emph{standard pairing} of $G$ as
$$\left<(a_1,\dots,a_m),(b_1,\dots,b_m)\right> = \zeta_{n_1}^{a_1b_1}\zeta_{n_2}^{a_2b_2}\dots\zeta_{n_m}^{a_mb_m}$$
where $\zeta_n = e^{2\pi i/n}$.

An important tool in the study of formal duality is the reduction to so called primitive sets. A primitive set is defined as follows.

\begin{definition}\label{def:primitive}
A set $S\subset G$ is called \emph{primitive} if \textbf{none} of the following holds:
\begin{enumerate}
\item $S \subset v\cdot H$ for some proper subgroup $H$ of $G$,
\item $S$ is a union of cosets with respect to a non-trivial subgroup of $G$.
\end{enumerate}
\end{definition}

The characterization of formally dual sets reduces to the study of primitive formally dual sets by the following result:

\begin{theorem}[{\cite[Lemma 4.2]{cohn2014formal}, \cite[Theorem 3.15]{schuler2019phd}}]\label{thm:reduce_to_primitive}
Let $S\subset G$ and $T\subset\hat G$ form a formally dual
pair. Then the following holds:
\begin{enumerate}
\item The set $S$ is contained in a proper coset $a\cdot H$ of $G$ if and only if $T$ is
a union of cosets with respect to $H^\perp\coloneqq \{\chi\in\hat G \ : \ \chi(h) = 1 \text{ for all } h\in H\}$ and vice versa.
\item If $S\subset a\cdot H$ then $S$ under the canonical map $a\cdot H \rightarrow H$ and $T$ under the natural reduction map $\hat G\rightarrow \hat G/H^\perp\simeq \hat H$ also form a formally dual pair.
\end{enumerate}
\end{theorem}

An equivalent definition of formal duality using the so called \emph{even sets} has been introduced in \cite[Section 4]{li2018abelian}. For this purpose, we consider the group algebra $\Q G$. More information about group algebras can be found in \cite[page 104]{serge2002algebra}.
Here, for a finite abelian group $(G, \cdot)$, the \emph{group algebra} $\Q G$ is the set of formal sums $A = \sum_{g\in G} a_g g$ with coefficients $a_g\in\Q$.
We define addition and multiplication in $\Q G$ by
$$\left( \sum_{g\in G} a_g g\right) + \left(\sum_{g\in G} b_g g\right) = \sum_{g\in G} (a_g + b_g) g, $$
$$\left( \sum_{g\in G} a_g g\right) \cdot \left(\sum_{g\in G} b_g g\right) = \sum_{g\in G}\sum_{h\in G} (a_g a_h) g\cdot h = \sum_{g\in G} \left(\sum_{h\in G} a_h a_{h^{-1}g}\right)g.$$
With this addition and multiplication $\Q G$ is indeed an algebra.
Also, for $A = \sum_{g\in G} a_g g\in\Q G$ we use $A^{(-1)} \coloneqq \sum_{g\in G} a_g g^{-1}$.

Furthermore, we abuse notation as $S = \sum_{g\in S} g\in\Q G$  as well as $S^{(-1)} \coloneqq \sum_{g\in S} g^{-1}\in\Q G$ for $S\subset G$. Note that this abuse of notation can also be applied to subgroups of $G$. Also, note that 
$$SS^{(-1)} = \sum_{g\in S}\sum_{h\in S}gh^{-1} = \sum_{g\in G} \nu_S(g)g.$$
Then, an \emph{even set} $S$ is a set such that $SS^{(-1)}$ can be expressed as a linear combination of subgroups in the group algebra using the mentioned abuse of notation. For more details about even sets kindly refer to \cite[Section 4]{li2018abelian} or \cite[Chapter 4]{schuler2019phd}. We recall the main theorem of this perspective:

\begin{theorem}[{\cite[Theorem 4.9]{li2018abelian}}]\label{thm:even_set_main}
Two sets $S\subset G$ and $T\subset \hat G$ form a formally dual pair if and only if there are coefficients $\lambda_H\in\Q$ such that
\begin{enumerate}
\item $SS^{(-1)} = \sum_{H\leq G} \lambda_H H$ and
\item $TT^{(-1)} = \sum_{H\leq G} \left(\frac{\lambda_H |G| \cdot |H|}{|S|^3}\right) H^\perp$.
\end{enumerate}
Especially $S$ and $T$ are even sets.
\end{theorem}

The following example states previously known formally self dual sets:

\begin{example}\label{ex:non-prim-fsd}
The following sets are formally self dual sets.
\begin{enumerate}
\item The trivial example $S = \{1\} \subset G = \{1\}$,
\item $\TITO \coloneqq \{0,1\}\subset \Z_4$ under the standard pairing (\cite[Section 3.1]{cohn2014formal}),
\item any $(n,n,n,1)$-relative difference set with forbidden subgroup $N$ such that $\Delta(N)=N^\perp$ under $\Delta$ (\cite[Theorem 3.7]{li2018abelian}),
\item generalized relative difference sets in $\Z_{p^k}^{2s}$ under a specific isomorphism (\cite[Section 3.3]{li2018abelian}),
\item The set $\{n\cdot k \ : \ k\in \Z_{n^2}\}$ under the standard pairing \cite[Theorem 6.1]{xia2016classification},\label{it:xia_ex_one}
\item For any prime $p$ and any $\alpha$ with $\alpha^2 \equiv -1 \mod p$ the set 
$$\{ (k,k\cdot \alpha) \ : \ k\in \Z_p\}$$
is formally self dual in $(\Z_p)^2$ under the standard pairing (by the choice of $\alpha$, we might only consider $p$ with $p\equiv 1 \mod 4$) \cite[Theorem 6.3]{xia2016classification}.\label{it:xia_ex_two}
\end{enumerate}
\end{example}

There is another way known to construct primitive formally dual sets of the same size using skew Hadamard difference sets:

\begin{definition}
Let $G$ be a finite abelian group. A set $D$ is called a \emph{skew Hadamard difference set} if
\begin{enumerate}
\item $DD^{(-1)} = \lambda G + k$ for some $\lambda,k\in\Z$,
\item $G = D + D^{(-1)} + \{1\}$ ($D$, $D^{(-1)}$ and $\{1\}$ are a partition of $G$).
\end{enumerate}
\end{definition}

\begin{theorem}[{\cite[Theorem 3.20]{li2018abelian}}]\label{thm:sHDS_construction}
Let $p^m\equiv 3 \pmod 4$ be a prime power and $D$ be a skew Hadamard difference set in $\Z_p^m$ and $\Delta:\Z_p^m\rightarrow \hat{\Z_p^m}$ be a group isomorphism.
Let
$$D^* = \left\{a \in \Z_p^m \ : \ \left<a,D\right>_\Delta = \frac{-1 + i\sqrt{p^m}}{2}\right\}$$
be the dual set of $D$ (see \cite{li2018abelian} and \cite{weng2009skewHadamardDifferneceSets} for details).
 Furthermore, let $\alpha$ and $\beta$ be nonzero elements of $\Z_p$ with $\alpha\neq\beta$ and $G = \Z_p^m\times\Z_p^m$. Then
\begin{align*}
S &= \{(0,0)\} \cup \{(x,\alpha x) \ : \ x\in D\} \cup \{(x,\beta x) \ : \ x\in D^{(-1)}\}\\ \text{ and}\\
T &= \{(0,0)\} \cup \{(\frac{\alpha}{\alpha-\beta}x,\frac{1}{\beta - \alpha}x) \ : \ x\in D^*\} \cup \{(\frac{\beta}{\alpha-\beta}x, \frac{1}{\beta - \alpha}x) \ : \ x\in D^{*(-1)}\}
\end{align*}
form a formally dual pair under the isomorphism $(x,y)\mapsto (\Delta(x),\Delta(y))$.
\end{theorem}	

\begin{remark}
Theorem \ref{thm:sHDS_construction} can be easily obtained by {\cite[Theorem 3.20]{li2018abelian}}, even though it is a slight generalization of this result.
\end{remark}

In Theorem \ref{thm:paley_is_self_dual} we will prove that, at least in a special case, Theorem \ref{thm:sHDS_construction} produces a formally self dual set.
Therefore, consider the following example of a skew Hadamard difference set:

\begin{example}[{\cite[Example 2.5]{xiang2005algebraic_design_theory}}]
Let $p^m\equiv 3\pmod 4$ be a prime power. Then the set
$$D = \{x^2 \ : \ x\in \F_{p^m}\setminus \{0\}\}$$
is a skew Hadamard difference set in the additive group of $\F_{p^m}$ called the \emph{Paley difference set}.
\end{example}

\begin{theorem}\label{thm:paley_is_self_dual}
Let $D$ be the Paley difference set in $\Z_p^m$ and $\alpha,\beta\in\Z_p$ with $\alpha\neq\beta$.
Then
$$S = \{(0,0)\} \cup \{(x,\alpha x) \ : \ x\in D\} \cup \{(x,\beta x) \ : \ x\in D^{(-1)}\}$$
is a formally self dual set.
\end{theorem}
\begin{proof}
Let $\Tr$ be the field trace of $\F_{p^m}$. We consider the isomorphism $\Delta_{\Tr}$ with the respective pairing
$\left<a,b\right>_{\Delta_{\Tr}} = \zeta_p^{\Tr(a\cdot b)}$.

For $a\neq 0$ we have
\begin{equation}\label{eq:payley_trace}
\left<a,D\right>_{\Delta_{\Tr}} = \frac 12 \left(\sum_{x\in \F_{p^m}} \left<a,x^2\right>-1\right).
\end{equation}

It is known (see for example \cite[Theorems 5.15, 5.33]{niederreiter1996finite_fields}) that
$$\sum_{x\in \F_{p^m}} \left<a,x^2\right>_{\Delta_{\Tr}} = \begin{cases} \eta(a)(-1)^{m-1}\sqrt{p^m} \text{ if } p\equiv 1\pmod 4\\ \eta(a)(-1)^{m-1}i^m\sqrt{p^m} \text{ if } p\equiv 3\pmod 4\end{cases}$$

where $\eta$ is the quadratic character ($\eta(a) = 1$ if $a$ is square and $\eta(a) = -1$ if $a$ is nonsquare).

In the considered case, we have $p^m\equiv 3\pmod 4$ and therefore $p\equiv 3\pmod 4$ and $m\equiv 1 \pmod 2$.
Therefore
$$\sum_{x\in \F_{p^m}} \left<a,x^2\right>_{\Delta_{\Tr}} = \pm \eta(a)i\sqrt{p^m}.$$
By substituting this term in Equation \ref{eq:payley_trace} we get
$$\left<a,D\right>_{\Delta_{\Tr}} = \frac 12 \left(-1 \pm \eta(a)i\sqrt{p^m}\right).$$
This yields
\begin{equation}
D^* = \left\{a\in\F_{p^m} \ : \ \left<a,D\right>_{\Delta_{\Tr}} = \frac{-1 + i\sqrt{p^m}}{2}\right\} \in \{D, D^{(-1)}\}.
\end{equation}\label{eq:Dstar_under_trace}
By Theorem \ref{thm:sHDS_construction} we know that
\begin{align*}
S &= \{(0,0)\} \cup \{(x,\alpha x) \ : \ x\in D\} \cup \{(x,\beta x) \ : \ x\in D^{(-1)}\}\\ \text{ and}\\
T &= \{(0,0)\} \cup \left\{\left(\frac{\alpha}{\alpha-\beta}x,\frac{1}{\beta - \alpha}x\right) \ : \ x\in D^*\right\} \cup \left\{\left(\frac{\beta}{\alpha-\beta}x, \frac{1}{\beta - \alpha}x\right) \ : \ x\in D^{*(-1)}\right\}
\end{align*}
form a formally dual pair under $\Delta_{\Tr}\times\Delta_{\Tr}: (x,y)\mapsto(\Delta_{\Tr}(x), \Delta_{\Tr}(y))$.

Due to Equation \eqref{eq:Dstar_under_trace} this yields
$T = \pi(S)$
where 
$$\pi(x,y) = \begin{cases}\left(\frac{1}{\alpha-\beta}y, \frac{1}{\beta-\alpha}x\right) \text{ if } D^* = D\\
																\left(\frac{\alpha+\beta}{\alpha-\beta}x + \frac{1}{\beta-\alpha}y, \frac{1}{\beta-\alpha}x\right) \text{ if } D^* = D^{(-1)}\end{cases}$$
																
Therefore, $S$ is a formally self dual set under $\Delta_{\Tr}\times\Delta_{\Tr}\circ\pi$.															
\end{proof}

\begin{remark}
There are several other constructions of skew-Hadamard difference sets (see for example \cite{ding2015skew_hadamard_from_dickson}, \cite{ding2007skew_hadamard_from_Ree_Tits}, \cite{ding2006skew_hadamard_family}, \cite{feng2012cyclotomic_skew_hadamard}, \cite{weng2007pseudo_paley}).
It is unknown if these examples and Theorem \ref{thm:sHDS_construction} can be used to produce more formally self dual sets.
Another family of examples has been constructed in \cite{li2019constructions} (see also \cite{li2020direct_construction}). These examples are not formally self dual since the respective sets have unequal sizes. Also, it is unknown if there are any relative difference sets with forbidden subgroup $N$ such that $N$ and $N^\perp$ are not isomorphic. If such relative difference sets exist, it is unclear if they are formally self dual sets or even if they are formally dual sets.
\end{remark}

\section{Structural results about formally self dual sets}\label{sec:results}
It is easy to see that Examples \ref{ex:non-prim-fsd}.\ref{it:xia_ex_one} and \ref{ex:non-prim-fsd}.\ref{it:xia_ex_two} are not primitive. Thus, by Theorem \ref{thm:reduce_to_primitive} they can be reduced to a primitive formally dual set. However, it is not guaranteed that the resulting set is still a formally self dual set. The following proposition implies, that a result similar to Theorem \ref{thm:reduce_to_primitive} also holds for formally self dual sets, i.e., every non-primitive formally self dual set can be reduced to a primitive formally self dual set. A weaker version of this Proposition has been stated in the second authors thesis \cite[Proposition 3.22]{schuler2019phd}.

\begin{proposition}\label{prop:formal-self-duality}
Suppose $S\subset G$ is a non-primitive formally self dual set under an isomorphism $\Delta$. Then $S$ is contained in some proper coset of $H < G$. Furthermore, if $H$ is chosen as small as possible,
then $\tilde H \coloneqq \Delta^{-1} H^\perp \leq H$ and 
$$S' = \{v\tilde H \ : \ v\in S\} \subset H/\tilde H$$is formally self dual under the isomorphism $\Delta'$ given by
$\left<a \tilde H,b \tilde H\right>_{\Delta'} = \left<a,b\right>_\Delta.$
\end{proposition}
\begin{proof}
First suppose $S$ is not primitive but also not contained in a proper coset. Due to Definition \ref{def:primitive} we know that $S = \bigcup v\cdot L$ for some non-trivial $L\leq G$. Using Theorem~\ref{thm:reduce_to_primitive} we then know that $\Delta S$ is contained in a coset with respect to $L^\perp < \hat G$ and thus $S$ is contained in a coset with respect to $\Delta^{-1} L^\perp < G$ which is a contradiction.

Now assume that $S$ is contained in a coset with respect to $H<G$. By Theorem \ref{thm:reduce_to_primitive} we know that $\Delta S$ is invariant under translations by $H^\perp$ and thus $S$ is invariant under translations by $\tilde H$.

Therefore, using Theorem \ref{thm:reduce_to_primitive} again, we have $\Delta S \subset \tilde H^\perp$ and thus $S\subset H\cap \Delta^{-1}\tilde H^\perp$. Since $H$ is chosen as small as possible and $|H| = |\Delta^{-1}\tilde H^\perp|$ we have $H = \Delta^{-1}\tilde H^\perp$. In other words $\Delta H = \tilde H^\perp$.

This yields, that $\Delta'$ is well defined, since for all $h\in \tilde H$ and $b\in H$ we have 
\begin{enumerate}
\item $\left<h,b\right>_{\Delta} = 1$ (indeed $\Delta h\in \Delta \tilde H = \Delta \Delta^{-1} H^\perp = H^\perp$),
\item $\left<b,h\right>_\Delta = 1$ (indeed $\Delta b \in \Delta H = \tilde H^\perp$).
\end{enumerate}
Then we have for all $a\tilde H\in H/\tilde H$ that $\nu_{S'}(a\tilde H) = \nu_{S}(a) / |\tilde H|$ as well as $|S'| = |S|/|\tilde H|$ and 
$$\left<a,S\right>_\Delta = \sum_{v\in S} \left<a,v\right>_\Delta = |\tilde H|\cdot \sum_{v\tilde H\in S'} \left<a \tilde H,v\tilde H\right>_{\Delta'} = |\tilde H|\cdot \left <a\tilde H, S'\right>_{\Delta'}.$$
Altogether, we have
$$\frac{|S'|^2}{|S'|} \nu_{S'}(a\tilde H) = \frac 1{|\tilde H|^2} \frac{|S|^2}{|S|} \nu_S(a) = \frac 1 {|\tilde H|^2} \left| \left<a,S\right>_\Delta\right|^2 = \left|\left <a\tilde H, S'\right>_{\Delta'}\right|^2.$$
Thus, $S'$ is formally self dual under $\Delta'$.
\end{proof}

Due to this result, we only need to characterize the primitive formally self dual sets in order to characterize all formally self dual sets. In particular we have:

\begin{corollary}
By using Proposition \ref{prop:formal-self-duality} on Example \ref{ex:non-prim-fsd}.\ref{it:xia_ex_one} and Example \ref{ex:non-prim-fsd}.\ref{it:xia_ex_two} with $H = S$, the examples reduce to the trivial formally self dual set.
\end{corollary}

We continue by discussing formal self duality in perspective of even sets. Our main aim is to provide a canonical even set decomposition for formally self dual sets.
This is achieved in Proposition \ref{prop:even_set_fsd}. First we discuss some insights which are helpful in order to proof Proposition \ref{prop:even_set_fsd}.
Recall that an even set is a set $S\subset G$ such that $SS^{(-1)} = \sum_{H \leq G} \lambda_H H \in \Q G$ for some rational coefficients $\lambda_H$. For every isomorphism $\Delta: G\rightarrow \hat G$ we define the adjoint isomorphism by $\Delta_*: G\rightarrow \hat G$ such that $\left<x,y\right>_\Delta = \left<y,x\right>_{\Delta_*}$. In other words 
$$\chi(y) = [\Delta_* y](\Delta^{-1}\chi) = [\Delta y](\Delta_*^{-1}\chi)$$
for all $\chi\in\hat G$ and $y\in G$. This induces an automorphism $\sigma = \Delta^{-1}\Delta_*$ of $G$. Note, that if $\Delta$ corresponds to a symmetric bilinear form, for example the standard pairing, then $\Delta = \Delta_*$ and $\sigma = \id$. In the general case we can use the following to describe the relations among $H$, $H^\perp$, $\Delta$ and $\Delta_*$.

\begin{lemma}\label{lem:sigma_rules}
For some isomorphism $\Delta:G\rightarrow\hat G$ and subgroup $H\leq G$ we have:
\begin{enumerate}
\item $(\Delta_* H)^\perp = \Delta^{-1} H^\perp$ (see also \cite[Lemma 2.1]{schuler2019phd}),\label{it:Delta}
\item $\Delta_*^{-1}(\Delta^{-1} H^\perp)^\perp = H$, \label{it:Delta_diff}
\item $\Delta^{-1}(\Delta^{-1} H^\perp)^\perp = \sigma H$,
\item $\Delta^{-1}(\sigma H)^\perp = \sigma \Delta^{-1} H^\perp$.
\end{enumerate}
\end{lemma}
\begin{proof}
Let $h\in \Delta^{-1} H^\perp$. Then for all $\chi\in\Delta_* H$ we have
$$h(\chi) = \chi(h) = [\underbrace{\Delta h}_{\in H^\perp}](\underbrace{\Delta_*^{-1}\chi}_{\in H}) = 1.$$
Thus $h\in (\Delta_* H)^\perp$.

On the other hand let $g\in (\Delta_* H)^\perp$. Then for all $h\in H$ we have
$$[\Delta g](h) = [\Delta_* h](g) = g(\underbrace{\Delta_* h}_{\in \Delta_* H}) = 1.$$
Thus $\Delta g\in H^\perp$ or in other words $g\in\Delta^{-1}H^\perp$.
Altogether, we get assertion \eqref{it:Delta}.

Furthermore, $\Delta_*^{-1}(\Delta^{-1} H^\perp)^\perp = \Delta_*^{-1} \Delta_* H = H$ and
$$\Delta^{-1}(\Delta^{-1} H^\perp)^\perp = \Delta^{-1}\Delta_* H = \sigma H$$
as well as
$$\Delta^{-1}\underbrace{(\Delta^{-1}\Delta_* H)^\perp}_{=\Delta_*(\Delta_* H)^\perp} =  \sigma(\Delta_* H)^\perp = \sigma \Delta^{-1} H^\perp$$ yielding the other assertions.
\end{proof}

Furthermore, we have the following result about formal duality (see also \cite[Lemma 3.21]{schuler2019phd}):
\begin{lemma}\label{lem:Delta_to_Delta_ast}
Two sets $S\subset G$ and $T\subset G$ form a formally dual pair under isomorphism $\Delta$ if and only if $T$ and $S$ form a formally dual pair under $\Delta_*$. In particular, $S$ is formally self dual under $\Delta$ if and only if $S$ is formally self dual under $\Delta_*$.
\end{lemma}
\begin{proof}
Suppose $S$ and $T$ form a formally dual pair under $\Delta$. 
Observe that \mbox{$\nu_{\Delta T}(\Delta g) = \nu_T(g)$} and
$$[\Delta g](S) = \sum_{x\in S} [\Delta g](x) = \sum_{x\in S} [\Delta_* x](g) = \sum_{x\in S} g(\Delta_* x) = g(\Delta_* S).$$
Then
$$\frac{|\Delta_* S|^2}{|T|} \nu_T(g) = \frac{|S|^2}{|\Delta T|}\nu_{\Delta T}(\Delta g) = |[\Delta g](S)|^2 = |g(\Delta_* S)|^2.$$
With $[\Delta_*]_* = \Delta$ the assertion follows by Definition \ref{def:fd_under_isomorphism}.
\end{proof}

With these results we can provide a canonical even set decomposition of formally self dual sets:

\begin{proposition}\label{prop:even_set_fsd}
Let $S\subset G$ and $\Delta: G\rightarrow \hat G$ be an isomorphism.

The following are equivalent:
\begin{enumerate}
\item $S$ is a formally self dual set under $\Delta$,
\item $SS^{(-1)} = \sum_{H\leq G} \lambda_H H$ with $\lambda_{\sigma H} = \lambda_H$ and $\lambda_{\Delta^{-1} H^\perp} = \frac{|H|}{|S|}\lambda_H$.
\end{enumerate}
\end{proposition}
\begin{proof}
Suppose $S$ is a formally self dual set under $\Delta$. By Theorem \ref{thm:even_set_main} 
$$SS^{(-1)} = \sum_{H\leq G} \mu_H H$$
for suitable $\mu_H\in\Q$.

First we show by induction in $i$ that
\begin{equation}SS^{(-1)} = \sum_{H\leq G} \mu_{\sigma^i H} H = \sum_{H\leq G} \frac{|S|}{|H|} \mu_{\sigma^{i}\Delta^{-1} H^\perp} H \text{ for all }i\in\N.\label{eq:multi_sigma}\end{equation}
If $i=0$ the first equation follows from the above assumption $SS^{(-1)} = \sum_{H\leq G} \mu_H H$. For the second equation we use  Theorem \ref{thm:even_set_main}, Lemma~\ref{lem:Delta_to_Delta_ast} and $|S|^2 = |G|$ to get
$$SS^{(-1)} = \sum_{H\leq G} \underbrace{\frac{|G|}{|S|^3}|H|}_{=\frac{|H|}{|S|}} \mu_H \Delta_*^{-1} H^\perp.$$
By changing the order of summation by substituting $H$ by $\Delta^{-1} H^\perp$ and using Lemma \ref{lem:sigma_rules} we get
$$SS^{(-1)} = \sum_{H\leq G} \underbrace{\frac{|\Delta^{-1} H^\perp|}{|S|}}_{\frac{|S|}{|H|}}\mu_{\Delta^{-1} H^\perp}\underbrace{\Delta_*^{-1} (\Delta^{-1} H^\perp)^\perp}_{= H}.$$
The proof of induction can be obtained in a similar fashion using Lemma \ref{lem:sigma_rules} as
\begin{align}
SS^{(-1)} 
&= \sum_{H\leq G} \frac{|S|}{|H|} \mu_{\sigma^i\Delta^{-1}H^\perp} H 
&= \sum_{H\leq G}\mu_{\sigma^i\Delta^{-1}H^\perp}\Delta_*^{-1}H^\perp \\
&= \sum_{H\leq G} \underbrace{\mu_{\sigma^i\Delta^{-1}(\Delta^{-1} H^\perp)^\perp}}_{= \mu_{\sigma^{i+1} H}} H&\label{eq:multi_sigma_first_eq}\\
&= \sum_{H\leq G} \frac{|H|}{|S|} \mu_{\sigma^{i+1}H}\Delta_*^{-1}H^\perp
&= \sum_{H\leq G} \frac{|S|}{|H|} \mu_{\sigma^{i+1}\Delta^{-1} H^\perp} H\label{eq:multi_sigma_sec_eq}.
\end{align}
Equation \eqref{eq:multi_sigma} then follows by the lines \eqref{eq:multi_sigma_first_eq} and \eqref{eq:multi_sigma_sec_eq}.

Now set $s = \ord{\sigma}$ and by taking the average over all equations of the form \eqref{eq:multi_sigma} we get
$SS^{(-1)} = \sum_{H\leq G} \lambda_H H$ with
$$\lambda_H \coloneqq \frac{1}{2s}\sum_{i=0}^{s-1}\left(\mu_{\sigma^{i} H} + \frac{|S|}{|H|} \mu_{\sigma^{i}\Delta^{-1}H^\perp}\right).$$
Observe that by Lemma \ref{lem:sigma_rules} and $\sigma^s = \sigma^0$:
$$\lambda_{\sigma H} = \frac{1}{2s} \sum_{i=0}^{s-1} \left(\mu_{\sigma^{i+1} H} + \frac{|S|}{|H|} \underbrace{\mu_{\sigma^{i}\Delta^{-1}(\sigma H)^\perp}}_{=\mu_{\sigma^{i+1}\Delta^{-1}H^\perp}}\right) = \lambda_H$$
and
$$\lambda_{\Delta^{-1} H^\perp} = \frac{1}{2s} \sum_{i=0}^{s-1} \left(\mu_{\sigma^i \Delta^{-1} H^\perp} + \frac{|H|}{|S|}\underbrace{\mu_{\sigma^i\Delta^{-1}(\Delta^{-1} H^\perp)^\perp}}_{=\mu_{\sigma^{i+1}H}}\right) = \frac{|H|}{|S|}\lambda_H.$$
On the other hand, if $SS^{(-1)} = \sum_{H\leq G} \lambda_H H$ with $\lambda_H$ as asserted, then
$$SS^{-1} = \sum_{H\leq G} \frac{|H|}{|S|} \lambda_{\Delta^{-1} H^\perp} H = \sum_{H\leq G} \frac{|S|}{|H|} \lambda_H \Delta_*^{-1} H^\perp.$$
By Theorem \ref{thm:even_set_main} $S$ is formally self dual under $\Delta_*$ and thus by Lemma \ref{lem:Delta_to_Delta_ast} $S$ is formally self dual under $\Delta$.
\end{proof}

The following result discusses the relations between the canonical decomposition given in Proposition \ref{prop:even_set_fsd} and an arbitrary even set decomposition.

\begin{corollary}\label{cor:zero_sum}
Let $S$ be an even set and $SS^{(-1)} = \sum_{H\leq G} \mu_H H$. Furthermore, define $\lambda_H$ similar as in the proof of Proposition \ref{prop:even_set_fsd}. Then $S$ is formally self dual under $\Delta$ if and only if
$$\sum_{H\leq G} \left(\mu_H - \lambda_H\right) H = 0.$$
\end{corollary}
\begin{proof}
If $S$ is formally self dual, then $SS^{-1} = \sum_{H\leq G} \mu_H H = \sum_{H\leq G} \lambda_H H$ following the proof of Proposition \ref{prop:even_set_fsd}. The assertion easily follows.
On the other hand, if $\sum_{H\leq G} \left(\mu_H - \lambda_H\right) H = 0$ then
$$SS^{-1} = \sum_{H\leq G} \mu_H H + \sum_{H\leq G} \left(\lambda_H - \mu_H\right) H = \sum_{H\leq G} \lambda_H H$$
and $S$ is formally self dual under $\Delta$ by Proposition \ref{prop:even_set_fsd}.
\end{proof}

\section{New examples}\label{sec:examples}

In this section we discuss two examples of formally self dual sets in groups of order $64$ found via computer search.

\begin{example}\label{ex:sporadic_examples}
The following sets are primitive formally self dual sets under the standard pairing:
$$\{(0,0,0),(0,0,1),(0,0,2),(0,0,5),(0,1,0),(0,3,0),(1,0,0),(1,2,6)\}\subset\Z_2\times\Z_4\times\Z_8,$$
\begin{align*}
\{(0,0,0,0),(0,0,0,1),(0,0,0,2),&(0,0,0,5),(0,0,1,0),\\&(0,1,0,0),(1,0,0,0),(1,1,1,6)\}\subset\Z_2^3\times\Z_8
\end{align*}
\end{example}

The correctness of these examples can easily be verified by a computer algebra system. However, the two examples have a common structure, which possibly contains important information.
\begin{remark}\label{rem:generalization_attempt}
Let $S$ be one of the examples in Example \ref{ex:sporadic_examples}. Then $G = H_0\times H_1\times H_2$ ($H_0 = \Z_2$, $H_2 = \Z_8$, $H_1 = \Z_4$ and $H_1 = \Z_2^2$ respectively).
Futhermore, there is a set $S_0\subset H_0$, an element $e$ ($(1,2,6)$ and $(1,1,1,6)$ respectively) as well as $H_1' < H_1$ such that
\begin{enumerate}
\item $S = S_0 + H_1 - H_1' + H_2 + e$,
\item $[H_1':H_1] = 2$,
\item $2e\in H_0$ (written additively),
\item $4e = 0$,
\item $S_0 + S_0^{-1} = H_0 - \left<2e\right>$,
\item $S_0 + e$ is an even set.
\end{enumerate}
This representation of $S$ can be used to find an even set representation 
$$SS^{-1} = \sum_{H\leq G} \lambda_H H$$ (with many nonzero coefficients). Then Corollary \ref{cor:zero_sum} can be used to find an equation $0 = \sum_{H\leq G} \mu_H H$ which, together
with the stated properties (possibly with the addition of a few more) could yield more examples of formal self duality. However, this is not understood in detail yet.
\end{remark}

\section{Formally self dual sets from vectorial Boolean functions}\label{sec:boolean_functions}

In this section, we investigate formally self dual sets of the form $\{(x,F(x))\colon x \in \F_{2^n}\}~\subset~\F_{2^n}^2$, where $\F_{2^n}$ denotes as usual the field with $2^n$ elements and $F:\F_{2^n}\rightarrow\F_{2^n}$ is a function (such functions are also called \emph{vectorial Boolean functions}). Recall that the additive group of $\F_{2^n}$ is isomorphic to $\Z_2^n$. We define the \emph{absolute trace mapping} of $\F_{2^n}$ as $\Tr \colon \F_{2^n} \rightarrow \F_2$ via $x \mapsto x+x^2+\dots x^{2^{n-1}}$. 

Our investigation is motivated by the following new sporadic example that has been found by Shuxing Li~\cite{li2020private} who gave his permission to publish it here.

\begin{example}\label{ex:power_map}
The set $S \coloneqq \{(x,x^3) \ : \ x\in\F_8\}$ is formally self dual under the trace pairing $\left<(x,y),(a,b)\right> =(-1)^{\Tr (ax + by)}$.
\end{example}

The correctness of this example can easily be verified by a computer algebra system. Of course, the natural question is whether this specific example generalizes to larger fields and if similar formally self dual sets of the form $\{(x,F(x))\colon x \in \F_{2^n}\} \subset \F_{2^n}^2$ under the trace pairing $\left<(x,y),(a,b)\right> =(-1)^{\Tr (ax + by)}$ exist. Let us start with some more general results.
\begin{definition}
	Let $F \colon \F_{2^n} \rightarrow \F_{2^n}$ be a vectorial Boolean function. We define 
	\begin{equation*}
		\delta_F(a,b)=|\{x \in \F_{2^n} \colon F(x+a)+F(x)=b\}|
	\end{equation*}
	for all $a,b \in \F_{2^n}$. Functions where $\delta_F(a,b) \leq 2$ for all $a \in \F_{2^n}^*$ and all $b \in \F_{2^n}$ are called almost perfect nonlinear (APN).
	
	Moreover, we define the Walsh transform of $F$ by
		\begin{equation*}
		W_F(a,b)=\sum_{x \in \F_{2^n}}(-1)^{\Tr(bF(x)+ax)}
	\end{equation*}
	for all $a,b \in \F_{2^n}$.
\end{definition}
Both the Walsh transform and the differential uniformity of vectorial Boolean functions have been the subject of much research because of their immense importance to symmetric cryptography: If $F$ is used as an S-box of a block cipher, the Walsh transform immediately determines the resistance of $F$ to linear attacks and $\max_{a \in \F_{2^n}^*,b\in\F_{2^n}}\delta_F(a,b)$ the resistance to differential attacks. For an overview on vectorial Boolean functions and their cryptographic properties, we refer the reader to~\cite{carletbook}. Note that $\delta_F(a,b)$ is always even. Indeed, if $x$ is a solution of $F(x+a)+F(x)=b$ then so is $x+a$. We can now characterize formally self dual sets of the structure in Example~\ref{ex:power_map} via the values of $\delta_F$ and the Walsh transform.

\begin{theorem}\label{thm:umformulierung}
	The set $S=\{(x,F(x))\colon x \in \F_{2^n}\}\subset \F_{2^n}^2$ with $F \colon \F_{2^n} \rightarrow \F_{2^n}$ is formally self dual under the pairing $\langle (x,y),(a,b) \rangle = (-1)^{\Tr(ax+by)}$ if and only if
	\begin{equation}
		2^n\delta_F(a,b)=(W_F(a,b))^2
	\label{eq:selfdual}
	\end{equation}
	holds for all $a,b \in \F_2^n$.
\end{theorem}
\begin{proof}
	We have $v_S(a,b)=|\{(x,y) \in \F_{2^n}^2 \colon x+y=a, F(x)+F(y)=b\}|$. By substituting $y=x+a$ in the second equation, we get 
	$$v_S(a,b)=|\{x \in \F_{2^n} \colon F(x)+F(x+a)=b\}|=\delta_F(a,b).$$
	
	Further, we have $\langle (a,b),S \rangle = \sum_{x \in \F_{2^n}}(-1)^{\Tr(bF(x)+ax)}=W_F(a,b)$. The necessary condition for formally self dual sets now becomes $2^n\delta_F(a,b)=(W_F(a,b))^2$ as stated in the theorem.
\end{proof}

\begin{remark}\label{rem:trivial}
	Of course, it is easy to recognize a trivial formally self dual set: Choosing $F(x)=x$ the set $S=\{(x,x)\colon x \in \F_{2^n}\}\subset \F_{2^n}^2$ is clearly formally self dual.
\end{remark}

\begin{corollary}
If $S=\{(x,F(x))\colon x \in \F_{2^n}\}\subset \F_{2^n}^2$ with $F \colon \F_{2^n} \rightarrow \F_{2^n}$ is a formally self dual set, then $F(x)$ is bijective.
\end{corollary} 
\begin{proof}
	Assume $S$ is formally self-dual. We have $\delta_F(0,b)=0$ for all $b \neq 0$. By Equation~\eqref{eq:selfdual}, then $W_F(0,b)=0$ for all $b \neq 0$, i.e.
	\begin{equation*}
		W_F(0,b)=\sum_{x \in \F_{2^n}}(-1)^{\Tr(bF(x))}=0.
	\end{equation*}
	It is well known that this holds if and only if $F(x)$ is a bijection (see e.g. \cite[Theorem 7.7.]{niederreiter1996finite_fields}).
\end{proof}
\begin{remark} \label{rem:bijective}
For bijective funtions, the case $a=0$ or $b=0$ in Equation \eqref{eq:selfdual} holds trivially. The case $a=0$, $b \neq 0$ was used in the corollary above, and if $b=0$ and $a \neq 0$ then $\delta_F(a,0)=W_F(a,0)=0$. If $a=b=0$ then $\delta_F(0,0) = W_F(0,0)=2^n$. 
\end{remark}
Note that this immediately shows us that the set $\{(x,x^3) \colon x \in \F_{2^n}\}$ cannot be formally self dual if $n$ is even as the mapping $x \mapsto x^3$ is bijective on $\F_{2^n}$ if and only if $n$ is odd.  

\begin{corollary}
	If  $S=\{(x,F(x))\colon x \in \F_{2^n}\}\subset \F_{2^n}^2$ with $\F \colon \F_{2^n} \rightarrow \F_{2^n}$ is a formally self dual set, then $W_F(a,b)$ is divisible by $2^{\lceil \frac{n+1}{2}\rceil}$ for all $a,b \in \F_{2^n}$.
\end{corollary}
\begin{proof}
		By Equation \eqref{eq:selfdual} we have that $2^n\delta_F(a,b)=(W_F(a,b))^2$ for all $a,b \in \F_{2^n}$. The result follows since $\delta_F(a,b)$ is always even.
\end{proof}

Generally, determining the precise values for $\delta_F(a,b)$ as well as for the Walsh transform using theoretic means is quite difficult. By Theorem~\ref{thm:umformulierung} the number of different values of $|W_F(a,b)|$ has to be the same as the number of different values of $\delta_F(a,b)$. For an arbitrary function $F$ this is not expected. There is however one special class of functions for which this is always satisfied.

\begin{definition}
	Let $F \colon \F_{2^n} \rightarrow \F_{2^n}$ be a function. We call $F$ almost bent (AB) if $W_F(a,b) \in \{0,\pm 2^{\frac{n+1}{2}}\}$ for all $a \in \F_{2^n}$ and $b \in \F_{2^n}^*$
\end{definition}
Of course, AB functions exist only if $n$ is odd since $W_F(a,b)$ is always an integer. It is well known that all AB functions are also APN, i.e. $\delta_F(a,b) \in \{0,2\}$ for all $a \in \F_{2^n}^*$ and $b \in \F_{2^n}$ \cite{CCZ1998}. 

\begin{theorem}\label{thm:ab}
	Let $n$ be odd and $F \colon \F_{2^n} \rightarrow \F_{2^n}$ be a bijective AB function.
	The set $S=\{(x,F(x))\colon x \in \F_{2^n}\}$ is formally self dual under the pairing $\langle (x,y),(a,b) \rangle = (-1)^{\Tr(ax+by)}$if and only if 
	\begin{equation} \label{eq:ab}
		\delta_F(a,b)=0 \iff W_F(a,b)=0.
	\end{equation}
	for all $a,b \in \F_{2^n}$.
\end{theorem}
\begin{proof}
	By Equation \eqref{eq:selfdual} we need to show that 
	\begin{equation}
		2^n\delta_F(a,b)=(W_F(a,b))^2
	\label{eq:bedingung}
	\end{equation}
	for all $a,b \in \F_{2^n}$. Since $F$ is bijective, the equation holds if $a=0$ or $b=0$ by the considerations above. Now assume that $\delta_F(a,b) \neq 0$. Since $F$ is APN this necessarily means $\delta_F(a,b)=2$. 
	Then by Equation \eqref{eq:bedingung} necessarily $W_F(a,b)=\pm2^{\frac{n+1}{2}}$. Since $F$ is almost bent, this is equivalent to $W_F(b,a)\neq 0$. 
\end{proof}

AB functions are rare and finding infinite families of AB functions remains a difficult research problem. A particularly well-studied group of AB functions are AB monomials. It is known that all AB monomials are bijections, so they are natural choices for our search for formally self dual pairs. A table of all known AB monomials is given in Table~\ref{t:AB}, for references we again refer to~\cite[Section 3.1.6.]{carletbook}. In the following, we refer to functions given by the exponents in Table~\ref{t:AB} as Gold-functions, Kasami functions etc. As we can see, the function $x \mapsto x^3$ that leads to Example~\ref{ex:power_map} can be identified as a Gold function.
\begin{table}[ht]
	\centering
	\begin{tabular}{ ||c|c|c|| } 
	 \hline
		& Exponent & Conditions \\
	 \hline  \hline
	 Gold & $2^r+1$ & $\gcd(r,n)=1$  \\
	 \hline
	 Kasami & $2^{2r}-2^r+1$ & $\gcd(r,n)=1$\\
	 \hline
	 Welch & $2^t+3$ & \\
	 \hline
	 Niho & $2^t-2^{\frac{t}{2}}-1$ & $t$ even \\
	 & $2^t-2^{\frac{3t+1}{2}}-1$ & $t$ odd\\
	 \hline
	\end{tabular}
	\caption{List of known AB exponents over $\F_{2^n}$ with $n=2t+1$ (up to equivalences)}
	\label{t:AB}
\end{table}

The precise values of $\delta_F(a,b)$ and of the Walsh transform of AB functions can only be determined easily in the case of Gold functions.

\begin{proposition} \label{prop:gold1}
	Let $F \colon \F_{2^n} \rightarrow \F_{2^n}$ be defined by $F(x) = x^{2^i+1}$ with $n$ odd and $\gcd(i,n)=1$. Then 
	\begin{equation*}
		|W_F(a,1)| = \begin{cases} 0, & \Tr(a)=0,\\
			 2^{\frac{n+1}{2}}, & else.
			\end{cases}
	\end{equation*}
	and 
		\begin{equation*}
		\delta_F(1,b) = \begin{cases} 0, & \Tr(b)=0,\\
			2, & else.
						\end{cases}
	\end{equation*}
\end{proposition}
\begin{proof}
	The first statement is a classical result by Gold \cite{Gold}. 
	
	For the second statement, expanding yields 
	\begin{equation*}
		(x+1)^{2^i+1}+x^{2^i+1}= x^{2^i}+x+1=b.
	\end{equation*}
	Note that $\Tr(x^{2^i}+x+1)=1$, so the equation has no solution if $\Tr(b)=0$. Since $\sum_{b \in \F_{2^n}}\delta_F(1,b)=2^n$ and $\delta_F(1,b) \in \{0,2\}$, we know that $\delta_F(1,b)=0$ for precisely $2^{n-1}$ choices of $b$, and the result follows.
\end{proof}

\begin{corollary}\label{cor:gold}
	Let $F \colon \F_{2^n} \rightarrow \F_{2^n}$ be defined by $F(x) = x^{2^i+1}$ with $n$ odd and $\gcd(i,n)=1$. Then 
	\begin{equation*}
		W_F(a,b)=0 \text{ if and only if } \Tr(ab^{-1/(2^i+1)})=0
	\end{equation*}
	and 
	\begin{equation*}
		\delta_F(a,b)=0 \text{ if and only if } \Tr(a^{-(2^i+1)}b)=0.
	\end{equation*}
	Here, $1/(2^i+1)$ denotes the inverse of $2^i+1$ modulo $2^n-1$.
\end{corollary}
\begin{proof}
Observe that
\begin{equation*}
		W_F(a,b)=W_F(ab^{-1/(2^i+1)},1)
	\end{equation*}
	by substituting $x \mapsto a^{-1/(2^i+1)}x$ (note that $2^i+1$ is invertible in $\Z_{2^n-1}$ since $\gcd(i,n)=1$), and the first equation follows from Proposition~\ref{prop:gold1}. 
	
	Further, we have 
	\begin{equation*}
		x^{2^i+1}+(x+a)^{2^i+1}=x^{2^i+1}+a^{2^i+1}\left(\frac{x}{a}+1\right)^{2^i+1}
	\end{equation*}
	and by substituting $x \mapsto ax$ we get that $\delta_F(a,b)$ is the number of solutions of 
	$$a^{2^i+1}(x^{2^i+1}+(x+1)^{2^i+1})=b,$$ or, equivalently, 
	$$x^{2^i+1}+(x+1)^{2^i+1}=ba^{-(2^i+1)},$$
	so $\delta_F(a,b) = \delta(1,ba^{-(2^i+1)})$ and the result follows again from Proposition~\ref{prop:gold1}. 
\end{proof}

 \begin{theorem}\label{thm:gold2}
Let $n$ odd and $F \colon \F_{2^n} \rightarrow \F_{2^n}$ be a Gold function defined by $x~\mapsto~x^{2^i+1}$ with $i<n$, $\gcd(i,n)=1$.
	The set $S=\{(x,F(x))\colon x \in \F_{2^n}\}$  is formally self dual under the pairing $\langle (x,y),(a,b) \rangle = (-1)^{\Tr(ax+by)}$ if and only if $n=3$. In this case, both $i=1,2$ yield formally self dual sets.
\end{theorem}
\begin{proof}
	Using Theorem~\ref{thm:ab} and Corollary~\ref{cor:gold}, we just need to check if 
	$$\Tr(ab^{-1/(2^i+1)})~=~\Tr(a^{-(2^i+1)}b)$$
	holds for all $a,b \in \F_{2^n}$. Define the polynomial 
	$$G_b(a)=\Tr(ab^{-1/(2^i+1)})+\Tr(a^{-(2^i+1)}b).$$
	$G_b(a)$ has to be identical to the zero function (as a function in $a$) to satisfy the condition. $\Tr(ab^{-1/(2^i+1)})$ is linear in $a$. On the other hand $\Tr(a^{-(2^i+1)}b)$ is only linear if $n-(2^i+1)$ is a power of two. Clearly, the binary weight of $n-(2^i+1)$ is $n-2$, so this is only possible for $n=3$. For $n=3$ it can easily be checked that both $i \in \{1,2\}$ satisfy the condition.
\end{proof}

Theorem~\ref{thm:gold2} shows that the Example~\ref{ex:power_map} is a sporadic example in the sense that the set $\{(x,x^3)\colon x \in \F_{2^n}\}\subset \F_{2^n}^2$ only yields a formally self dual set if $n=3$. Of course, it would be interesting to find other formally self dual sets of the form $\{(x,F(x))\colon x \in \F_{2^n}\}$. Natural candidates would be the other AB monomials listed in Table~\ref{t:AB}. However, computer searches in low dimensions do not yield any new examples, and a theoretic treatment of the condition in Theorem~\ref{thm:ab} is very difficult for the non-Gold AB monomials. 

We like to note that when we find a formally self dual set of the form $\{(x,F(x))\colon x \in \F_{2^n}\}$ as in Example~\ref{ex:power_map}, then it is possible to construct some other formally self dual sets of this form.

For a group $G$ and $\phi$ in its automorphism group $\Aut(G)$ we can define the \emph{adjoint} $\phi^*$ of $\phi$ such that $\langle \phi(x),y\rangle = \langle x,\phi^*(y) \rangle$ for all $x,y \in G$. Clearly, $\phi^* \in \Aut(G)$.

\begin{proposition}[{\cite[Proposition 2.16]{li2018abelian}}] \label{prop:equiv}
	Let $G$ be a group. Let $\phi \in \Aut(G)$ and $\phi^*$ be the adjoint of $\phi$. Suppose $S$ and $T$ form a formally dual pair in $G$. Then $\phi(S)$ and $(\phi^*)^{-1}(T)$ also form a formally dual pair in G.
\end{proposition}

We can state the following straightforward corollary of Proposition~\ref{prop:equiv}.
\begin{corollary}\label{cor:adjoint}
	Let $G$ be a group. Let $\phi \in \Aut(G)$ and $\phi^*$ be the adjoint of $\phi$. Suppose $S$ is a formally self dual pair in $G$ and $\phi = (\phi^*)^{-1}$. Then $\phi(S)$ is also a formally self dual pair in $G$.
\end{corollary}

\begin{remark}
	If we consider the group $G = \F_{2^n} \times \F_{2^n}$ with the component-wise operation $(a,b)+(c,d) = (a+c,b+d)$ and consider sets of the form $\{(x,F(x)) \colon x \in \F_{2^n}\}$, the construction in Corollary~\ref{cor:adjoint} is closely related to the concept of \emph{CCZ-equivalence} (sometimes also called graph-equivalence) of vectorial Boolean functions. We again refer the reader to~\cite{carletbook}.
\end{remark}

We will now give a construction how to use Corollary~\ref{cor:adjoint} to construct new formally self dual sets of the form $\{(x,F(x)) \colon x \in \F_{2^n}\}$ under the trace pairing $\left<(x,y),(a,b)\right> =(-1)^{\Tr (ax + by)}$ from known ones.

Let $G = \F_{2^n} \times \F_{2^n}$ be a group with the component-wise operation $(a,b)+(c,d) = (a+c,b+d)$ and $\phi \in \Aut(G)$. We can write $\phi=\begin{pmatrix}
	L_1 & L_2 \\
	L_3 & L_4
\end{pmatrix}$, where $L_1,L_2,L_3,L_4$ are $\F_2$-linear mappings from $\F_{2^n}$ to $\F_{2^n}$. Applying $\phi$ to a vector $(a,b) \in \F_{2^n}^2$ yields (similar to the usual matrix multiplication)
\[\phi(a,b) = \begin{pmatrix}
	L_1 & L_2 \\
	L_3 & L_4
\end{pmatrix} (a,b) = (L_1(a)+L_2(b),L_3(a)+L_4(b)).\] 

\begin{lemma} \label{lem:adjoint}
Let $\phi \in \Aut(G)$, written as 
$\phi=\begin{pmatrix}
	L_1 & L_2 \\
	L_3 & L_4
\end{pmatrix}$
where $L_1, L_2, L_3, L_4$ are $\F_2$-linear mappings from $\F_{2^n}$ to itself. Then
$\phi^*=\begin{pmatrix}
	L_1^* & L_3^* \\
	L_2^* & L_4^*
\end{pmatrix}$.
Here we denote by $L_i^*$ the adjoint of $L_i$ with respect to the trace bilinear form $(\cdot,\cdot) \colon \F_{2^n} \times \F_{2^n} \rightarrow \F_2$ defined by $(a,b) \mapsto \Tr(ab)$.
\end{lemma}
\begin{proof}
	Let $a_1,a_2,b_1,b_2 \in \F_{2^n}$. Then 
	\begin{align*}
		\langle \phi(a_1,b_1),(a_2,b_2)\rangle &= \Tr((L_1(a_1)+L_2(b_1))a_2)+\Tr((L_3(a_1)+L_4(b_1))b_2) \\
		&=\Tr(a_1(L_1^*(a_2)+L_3^*(b_2)))+\Tr(b_1(L_2^*(a_2)+L_4^*(b_2))),
	\end{align*}
	which directly implies the desired structure for $\phi^*$.
\end{proof}
Recall that every $\F_2$-linear mapping on $\F_{2^n}$ can be written as a polynomial of the form $L = \sum_{i=0}^{n-1}c_ix^{2^i} \in \F_{2^n}[x]$, and the adjoint of such a linear polynomial can be determined easily from the coefficients. Indeed, we have $L^*(x)=\sum_{i=0}^{n-1}c_i^{2^{n-i}}x^{2^{n-i}}$.

\begin{corollary}
	Let $F \colon \F_{2^n} \rightarrow \F_{2^n}$. If $S=\{(x,F(x))\colon x \in \F_{2^n}\}$ is a formally self dual set then so is $\{(x,F^{-1}(x))\colon x \in \F_{2^n}\}$, where $F^{-1}$ denotes the compositional inverse of $F$.
\end{corollary}
\begin{proof}
	We apply Corollary~\ref{cor:adjoint} with $\phi=\begin{pmatrix}
	L_1 & L_2 \\
	L_3 & L_4
\end{pmatrix}$ and $L_1=L_4=0$ and $L_2=L_3=x$. This choice implies $\phi=(\phi^*)^{-1}$. Indeed, by Lemma~\ref{lem:adjoint}, we have 
\[\phi^* = \begin{pmatrix}
	L_1^* & L_3^* \\
	L_2^* & L_4^*
\end{pmatrix} =  \begin{pmatrix}
	L_1 & L_2 \\
	L_3 & L_4
\end{pmatrix}=\phi,\] 
since $L_2(x)=L_3(x)=x$ and the zero-polynomial are self-adjoint. Observe that $\begin{pmatrix}
	0 & L \\
	L & 0
\end{pmatrix}$ for $L(x)=x$ is an involution, so $\phi=\phi^*=(\phi^*)^{-1}$ as claimed.
\end{proof}

\begin{proposition} \label{prop:affine}
Let $L_1,L_2$ be two bijective $\F_2$-linear mappings from $\F_{2^n}$ to itself. We set $F'=L_2\circ F \circ L_1$. Further assume $L_1=(L_1^{-1})^*$ and $L_2=(L_2^{-1})^*$.

If $\{(x,F(x))\colon x \in \F_{2^n}\}$ is a formally self dual set, then so is $\{(x,F'(x))\colon x \in \F_{2^n}\}$.
\end{proposition}
\begin{proof}
	Define $\phi=\begin{pmatrix}
	L_1^{-1} & 0 \\
	0 & L_2
\end{pmatrix}$.

Then $\phi(x,F(x))=(L_1^{-1}(x),L_2(F(x)))$. In particular, 
$$\{\phi(x,F(x)) \colon x\in \F_{2^n}\} = \{(x,F'(x)) \colon x\in \F_{2^n}\}.$$ The condition $\phi=(\phi^*)^{-1}$ can be easily verified, so the claimed result follows from Corollary~\ref{cor:adjoint}.
\end{proof}

In~\cite[Definition 2.17]{li2018abelian}, the formal dual pairs that can be constructed with Proposition~\ref{prop:equiv} (and, by extension, with Proposition~\ref{prop:affine}) are called \emph{equivalent} to the original one. Still, it is interesting to note that these equivalent formal self dual sets have the special structure of a graph, i.e. $\{(x,F(x))\colon x \in \F_{2^n}\}$. 

\begin{remark}
	With Remark~\ref{rem:trivial} and Proposition~\ref{prop:affine} we can create some more formally self dual sets. Indeed, if $L$ is a bijective linear function on $\F_{2^n}$ satisfying $L~=~(L^{-1})^*$, then $\{(x,L(x))\colon x \in \F_{2^n}\}$ is a formally self dual set. 
\end{remark}

We now give a criterion when the condition $L~=~(L^{-1})^*$ holds for some linear function $L \colon \F_{2^n} \rightarrow \F_{2^n}$.

\begin{proposition} \label{prop:condition}
	Let $L \colon \F_{2^n} \rightarrow \F_{2^n}$ be a linear function defined by $L(x)=\sum_{i=0}^{n-1}c_ix^{2^i}\in\F_{2^n}[x]$. The condition $L=(L^{-1})^*$ holds if and only if 
	\begin{align*}
		\sum_{i=0}^{n-1} c_i^{2^{n-i}}&=1 \\
		\sum_{i=0}^{n-1} (c_ic_{i+j})^{2^{n-i}} &=0
	\end{align*}
	for all $j \in \{1,\dots,n-1\}$.
\end{proposition}
\begin{proof}
		$L=(L^{-1})^*$ holds if and only if $L^*(L(x))=x$ for all $x \in \F_{2^n}$. A simple calculation yields:
		\begin{align*}
			L^*(L(x))&=\sum_{i=0}^{n-1}c_i^{2^{n-i}}(\sum_{j=0}^{n-1}c_jx^{2^j})^{2^{n-i}}\\
				&=\sum_{j=0}^{n-1}\sum_{i=0}^{n-1}(c_ic_j)^{2^{n-i}}x^{2^{j-i}} \\
				&=\sum_{j=0}^{n-1}\sum_{i=0}^{n-1} (c_ic_{i+j})^{2^{n-i}}x^{2^j}.
		\end{align*}
		The polynomial $L^*(L(x))$ has to be equal to the polynomial $x$, so all coefficients have to be $0$ except the coefficient for $j=0$, which has to be $1$.
\end{proof}
	
	This condition can be checked very easily with a computer. Note that the conditions always hold if $L(x)=x^{2^i}$. Moreover, if $L=(L^{-1})^*$, then also $L^*=((L^*)^{-1})^*$, so the condition always holds for $L$ and $L^*$ simultaneously. We now give an example of such a secondary construction from the formally self dual set $\{(x,x^3)\colon x \in \F_{8}\}\subset \F_8^2$.
	
	\begin{example}
		Let $a$ be a root of the irreducible polynomial $x^3+x+1 \in \F_{2}[x]$. We define $L(x)=(a+1)x+(a^2+a+1)x^2+(a^2+1)x^4 \in \F_8[x]$. It is easy to verify that the conditions in Proposition~\ref{prop:condition} hold, so we have $L=(L^{-1})^*$.  Further, let $F(x)=x^3$ be the cube function. By Proposition~\ref{prop:affine} we can use this polynomial to construct new formally self-dual sets in $\F_8^2$, e.g.
		\begin{itemize}
			\item $\{(x,L(F(x)))\colon x \in \F_{2^n}\}=\{(x,(a^2+a+1)x^6+(a^2+1)x^5+(a+1)x^3)\colon x \in \F_{2^3}\}$
			\item $\{(x,F(L(x)))\colon x \in \F_{2^n}\}=\{(x,(a^2+1)x^6+(a+1)x^5+(a^2+a)x^4+\newline(a^2+a+1)x^3+ax^2+a^2x)\colon x \in \F_{2^3}\}$
			\item $\{(x,L(F(L(x))))\colon x \in \F_{2^n}\}=\{(x,x^5+x^4+x)\colon x \in \F_{2^3}\}$.
		\end{itemize}
	\end{example}

\section{Open questions}\label{sec:open_questions}

We end this paper with some open questions which are valuable to study.

\begin{enumerate}
\item Are the formally dual sets obtained via Theorem \ref{thm:sHDS_construction} formally self dual sets if we use other skew Hadamard difference sets than the Paley difference sets?
\item If $S$ is a formally self dual set under isomorphism $\Delta$, can we describe \emph{all} isomorphisms under which $S$ is formally self dual?
\item Are there other formally self dual sets of the form $\{(x,F(x))\colon x \in \F_{2^n}\}$?
\item Are there other groups of order $64$ that contain primitive formally dual sets? The next smallest group where we neither can prove existence or non-existence of primitive formally dual sets would be of order $72$. See Table \ref{tab:order_64_results} in the appendix for the currently known results in groups of order $64$.
\end{enumerate}

\section*{Acknowledgement}
The authors like to thank Shuxing Li for introducing Example \ref{ex:power_map} and giving permission to publish it.
Furthermore, we are grateful to Frieder Ladisch for careful reading and much helpful advice. Also, we are grateful to Claude Carlet for 
bringing the concept of formal dual codes to our attention.
\bibliography{formal-self-duality}
\bibliographystyle{alpha}

\section*{Appendix}

\begin{table}[htp]
	\begin{tabular}{ ||c|c|c|c|| } 
	\hline
 $(|G|,|S|)$ & $G$ & prim. f.d.s. exists& reasoning\\
\hline
\hline
$(64,2)$ &  arbitrary & NO & {\cite[Corollary 5.11]{li2018abelian}}\\
\hline
 $(64,4)$ & $\Z_{2}^6$ & NO & {\cite[Proposition 5.9]{li2018abelian}}\\
 & $\Z_{4}\times\Z_{2}^4$ & NO & {\cite[Proposition 5.9]{li2018abelian}}\\
 & $\Z_{4}^2\times\Z_{2}^2$ & NO & {\cite[Proposition 5.9]{li2018abelian}}\\
 & $\Z_{4}^3$ & unknown & no result \\
 & $\Z_{8}\times\Z_{2}^3$ & NO & {\cite[Proposition 5.9]{li2018abelian}}\\
 & $\Z_{8}\times\Z_{4}\times\Z_{2}$ & unknown & no result \\
 & $\Z_{8}^2$ & unknown & no result \\
 & $\Z_{16}\times\Z_{2}^2$ & unknown & no result \\
 & $\Z_{16}\times\Z_{4}$ & unknown & no result \\
 & $\Z_{32}\times\Z_{2}$ & NO & {\cite[Proposition 5.10]{li2018abelian}}\\
 & $\Z_{64}$ &  NO & \cite[Section 4.2]{xia2016classification}\\
\hline
 $(64,8)$ & $\Z_{2}^6$ & YES & Example \ref{ex:power_map} \\
 & $\Z_{4}\times\Z_{2}^4$ & unknown & no result \\
 & $\Z_{4}^2\times\Z_{2}^2$ & unknown & no result \\
 & $\Z_{4}^3$ & YES &{ \cite[Proposition 3.2]{li2018abelian}}\\
 & $\Z_{8}\times\Z_{2}^3$ & YES & Example \ref{ex:sporadic_examples} \\
 & $\Z_{8}\times\Z_{4}\times\Z_{2}$ & YES & Example \ref{ex:sporadic_examples} \\
 & $\Z_{8}^2$ & unknown & no result \\
 & $\Z_{16}\times\Z_{2}^2$ & unknown & no result \\
 & $\Z_{16}\times\Z_{4}$ & unknown & no result \\
 & $\Z_{32}\times\Z_{2}$ & unknown & no result \\
 & $\Z_{64}$ & NO & \cite[Section 4.2]{xia2016classification}\\
\hline
	\end{tabular}
\caption{Known existence results about primitive formally dual subsets $S$ of groups $G$ of order $64$.}
\label{tab:order_64_results}
\end{table}

\end{document}